\documentclass[11.5pt]{amsart}
\usepackage{fullpage}
\usepackage{graphicx}
\usepackage{subfigure}
\usepackage{psfrag}
\usepackage{color}
\usepackage{setspace}
\newtheorem{theorem}{Theorem}
\newtheorem{remark}[theorem]{Remark}

\title{A nonlocal model for fluid-structure interaction with applications in hydraulic fracturing}

\author{D.~Z.~Turner}
\address{Correspondence to: Dr. Daniel Z. Turner, Department of Civil Engineering, 
  Stellenbosch University, Private Bag X1 Matieland 7602, South Africa. \emph{E-mail address:} {\ttfamily dzturner@sun.ac.za}, \emph{Phone: } +27 21 808 4434}
\date{\today}
\begin{document}

\begin{abstract}
Modeling important engineering problems related to flow-induced damage (in the context of hydraulic fracturing among others) depends critically on characterizing the interaction of porous media and interstitial fluid flow. This work presents a new formulation for incorporating the effects of pore pressure in a nonlocal representation of solid mechanics. The result is a framework for modeling fluid-structure interaction problems with the discontinuity capturing advantages of an integral based formulation. A number of numerical examples are used to show that the proposed formulation can be applied to measure the effect of leak-off during hydraulic fracturing as well as modeling consolidation of fluid saturated rock and surface subsidence caused by fluid extraction from a geologic reservoir. The formulation incorporates the effect of pore pressure in the constitutive description of the porous material in a way that is appropriate for nonlinear materials, easily implemented in existing codes, straightforward in its evaluation (no history dependence), and justifiable from first principles. A mixture theory approach is used (deviating only slightly where necessary) to motivate an alteration to the peridynamic pressure term based on the fluid pore pressure. The resulting formulation has a number of similarities to the effective stress principle developed by Terzaghi and Biot and close correspondence is shown between the proposed method and the classical effective stress principle.
\end{abstract}


\maketitle

\section{INTRODUCTION}
The use of hydraulic fracturing to increase yields from natural gas reservoirs has gained popularity recently due to changes in the economic climate, the desire to decrease dependence on foreign energy sources, and the fact that natural gas burns cleaner than coal when used to produce electricity. This renewed interest in hydraulic fracturing has presented a number of challenges associated with modeling the processes involved. In particular, the need exists to better model the evolution of damage in the reservoir to prevent unintentional penetration of underground water resources. This type of modeling also helps to establish injection parameters that will produce the best results and is useful in estimating yields from a particular well. Towards these modeling needs, this work presents a modified formulation for state-based peridynamics that includes the influence of fluid pore pressure in a form similar to the familiar principle of effective stress. Using the theory of interacting continua \cite{Bowen, Atkin1976, Bedford_Drumheller_IJES_1983_v21_p863, Rajagopal} we show that the principle of correspondence (the idea that the state-based peridynamic representation of the constitutive law is equivalent to the classical rendition) holds for the proposed formulation. This is an important result in that the proposed formulation preserves the material behavior of the porous media as represented in the classical theory and provides a way to experimentally obtain the material parameters for the proposed formulation in a manner consistent with the classical approach.

The development of the peridynamic theory in general is presented in \cite{Silling, Silling2, Silling3, Silling4} and has been applied to a number of engineering problems related to membranes, reinforced concrete, and carbon fiber nanotube reinforced composites. This work builds on the \emph{state-based} (as opposed to bond-based) version of the theory that is presented in \cite{Silling2}. The major advantages of the state-based approach include a material response that depends on collective quantities like volume change or shear angle and the ability to incorporate a constitutive model from the conventional theory of solid mechanics. The core contribution of this work is a modification of the force scalar state which introduces the effects of fluid pore pressure. We show, under certain conditions, that the proposed formulation is equivalent to the notion of effective stress as presented by Terzaghi and Biot \cite{Terzaghi, Biot, BiotWillis}, but when generalized can be used to incorporate other models for the interaction of the fluid and the porous media. In the present work we use an interaction model that is based on the drag produced by the fluid on the porous matrix, although the use of other models is equally feasible. To simplify the presentation, we assume that the fluid pressure has previously been evaluated numerically or analytically through some other means. In this way we avoid delving into the numerous complexities regarding a fully coupled analysis in which the fluid pressure is treated as an unknown. Although we make no mention of the stability of various coupling algorithms for this formulation, we intend to pursue this in a forthcoming work.

A general overview of the theory behind poroelasticity is presented in \cite{Detournay}. In subsurface geomechanical modeling poroelasticity, or the interaction of the pore fluid with the porous media plays an important role and has been the focus of several studies \cite{Dean, DeanSchmidt, Minkoff, Madhav, TurnerLDRD, MartinezLDRD}. In the present work we demonstrate the proposed method for a number of numerical examples motivated by problems found in \cite{deBoer, Diebels, Andersen, Dean}. In many cases analytic solutions are known for the particular problem. For the examples that do not have an easily derived analytic solution, the results of the proposed method are compared with a commercial finite element code in order to gauge accuracy.

An outline of the proceeding sections is as follows: We begin with the development of the modified state-based peridynamic theory and its numerical implementation. Then, we present the governing equations in their discrete form. To demonstrate the interesting features of the proposed method we then show its performance for several numerical examples. We then conclude with some general comments and discussion of the results.

\section{FORMULATION}
\begin{figure}[htb!]
	\centering
\includegraphics[scale=0.8]{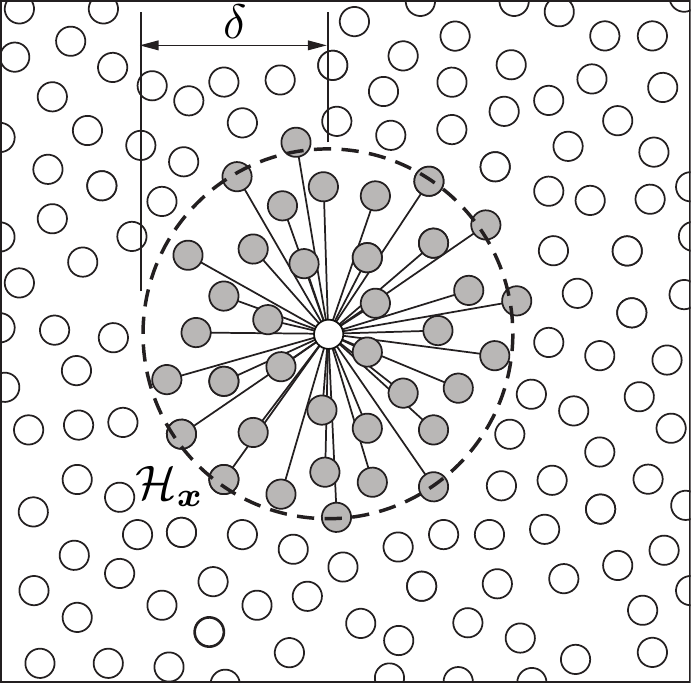}
	\caption{The spherical neighborhood, $\mathcal{H}_{\boldsymbol{x}}$, and horizon, $\delta$, of a particle at position $\boldsymbol{x}$ (the white particle at the center of the figure). The force interaction of this particle with the gray particles is considered non-negligible.}
	\label{fig:PeridynamicsNeighborhood}
\end{figure}

\begin{figure}[htb!]
	\centering
\includegraphics[scale=0.7]{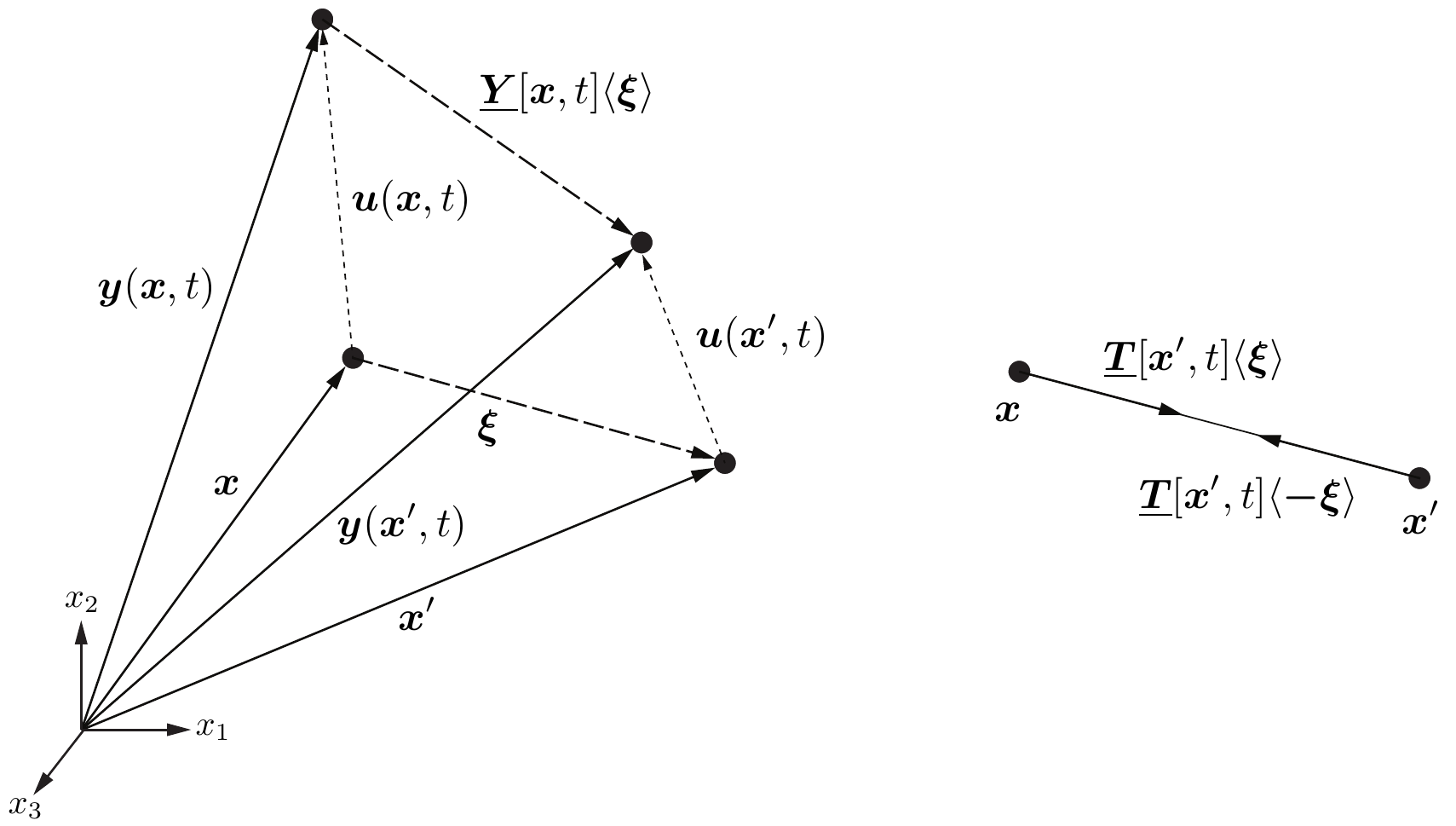}
	\caption{The geometry of deformation for peridynamics, (left) reference and deformed particle states (right) ordinary state based material response.}
	\label{fig:PeridynamicsGeometry}
\end{figure}

It will be helpful to begin with some basic definitions. The definition of a \emph{scalar state} and a \emph{vector state} are best understood in terms of their respective action upon a vector, $\boldsymbol{\xi}$. The scalar valued image of the vector, $\boldsymbol{\xi}$ under scalar state $\underline{s}$ is written: $\underline{s}\langle\boldsymbol{\xi}\rangle$. Likewise the vector valued image of the vector $\boldsymbol{\xi}$ under the vector sate $\underline{\boldsymbol{V}}$ is written $\underline{\boldsymbol{V}}\langle\boldsymbol{\xi}\rangle$. In general a state of order $m$ is a function $\underline{\boldsymbol{A}}\langle \cdot \rangle \; : \; \mathcal{H} \rightarrow \mathcal{L}_m$ where $\mathcal{H}$ is a spherical neighborhood centered at the origin in $\mathbb{R}^3$ and $\mathcal{L}_m$ is the set of all tensors of order $m$. The governing solid mechanics equations, as posed in the state-based peridynamic theory, make frequent use of the definitions above. In this work we restrict our attention to quasi-static analysis and neglect the inertial term (although the formulation does not preclude its inclusion). Consider a neighborhood, $\mathcal{H}_{\boldsymbol{x}}$, of a particle with position, $\boldsymbol{x}$ (see Figure~\ref{fig:PeridynamicsNeighborhood}). The state based equation of motion can be written in terms of a force vector state field, $\underline{\boldsymbol{T}}$ that defines the collective force relationship of the particle at $\boldsymbol{x}$ with each neighbor at position $\boldsymbol{x}'$.
\begin{align}
\int_{\mathcal{H}_x} \left\{ \underline{\boldsymbol{T}}[\boldsymbol{x},t]\langle \boldsymbol{x}' - \boldsymbol{x}\rangle - \underline{\boldsymbol{T}}[\boldsymbol{x}',t]\langle \boldsymbol{x} - \boldsymbol{x}'\rangle \right\} \mathrm{d}V_{\boldsymbol{x}'} + \boldsymbol{b}(\boldsymbol{x},t) = 0
\end{align}
The radius of the spherical neighborhood $\mathcal{H}_{\boldsymbol{x}}$, referred to as the \emph{horizon},  will be denoted $\delta$. The geometry of deformation is shown in Figure~\ref{fig:PeridynamicsGeometry}.
Regarding the integration in the equation above, in this formulation we assume that all functions are Riemann integrable. The relative position of two particles in the reference configuration is denoted as $\boldsymbol{\xi}$
\begin{align}
\boldsymbol{\xi} = \boldsymbol{x}' - \boldsymbol{x}
\end{align}
Let $\boldsymbol{y}(\boldsymbol{x},t)$ represent the the deformed position of the particle originally at $\boldsymbol{x}$
\begin{align}
\boldsymbol{y}(\boldsymbol{x},t) = \boldsymbol{x} + \boldsymbol{u}(\boldsymbol{x},t)
\end{align}
where $\boldsymbol{u}(\boldsymbol{x},t)$ is the displacement. The reference position and deformation vector state fields, $\underline{\boldsymbol{X}}$ and $\underline{\boldsymbol{Y}}$ respectively, are defined as
\begin{align}
\underline{\boldsymbol{X}}\langle \boldsymbol{\xi}\rangle &=  \boldsymbol{\xi}\\
\underline{\boldsymbol{Y}}[\boldsymbol{x},t]\langle \boldsymbol{\xi} \rangle = \boldsymbol{y}(\boldsymbol{x} &+ \boldsymbol{\xi},t) - \boldsymbol{y}(\boldsymbol{x},t)
\end{align}

\subsection{Linear elastic peridynamic solid}
For an ordinary material, there exists a scalar state, $\underline{t}$, acting in the direction of $\underline{\boldsymbol{Y}}$ such that the following relationship holds.
\begin{align}
\label{eq:ForceStateDef}
\underline{\boldsymbol{T}} = \underline{t} \;\underline{\boldsymbol{M}}
\end{align}
where $\underline{\boldsymbol{M}}(\underline{\boldsymbol{Y}})\langle \boldsymbol{x}' - \boldsymbol{x}\rangle$ is a unit vector pointing from the deformed position of $\boldsymbol{x}$ to the deformed position of $\boldsymbol{x}'$. The force scalar state for a linear elastic peridynamic material is given by
\begin{align}
\underline{t} = \frac{-3p}{m}\underline{\omega} \; \underline{x} + \beta \underline{\omega} \; \underline{e}^d
\end{align}
where $p = -k\theta$ is the peridynamic pressure, $\underline{x}$ is the reference position scalar state field, $\underline{e}^d$ is the deviatoric extension scalar state field, $\underline{\omega}$ is an influence function (which for this work we assume is of unit value), $m$ is the weighted volume, $\theta$ is the dilatation,  and $k$ and $\beta = 15\mu/m$ are material parameters. 
Via the principle of correspondence, it can be shown that $k=K^{(s)}$ represents the bulk modulus from the classical theory and $\mu = \mu^{(s)}$ represents the shear modulus.

The reference position scalar state field, $\underline{x}$, is simply the magnitude of the reference position vector state field. Likewise, let $\underline{y}$ represent the magnitude of $\underline{\boldsymbol{Y}}$. The deviatoric extension scalar state field, $\underline{e}^d$, stems from the extension scalar state field, $\underline{e} = \underline{y} - \underline{x}$, as follows
\begin{align}
\underline{e}^d = \underline{e} - \frac{\theta \underline{x}}{3}
\end{align}

\subsection{Effective stress principle}
Thus far we have made no deviation from the formulation presented in \cite{Silling2}) for an ordinary linear elastic peridynamic material. We now present a new peridynamic material, representative of a saturated porous medium effectively loaded by the fluid pressure in the pores. Since the force scalar state in the standard linear elastic peridynamic material above has been split into a dilatational or volume changing term and a deviatoric term, the extension to a peridynamic material with effective stress is straightforward. The influence of the fluid pore pressure can be included by modifying the definition of the peridynamic pressure as follows
\begin{align}
\label{eq:ModifiedConstit}
p = -k\theta + \gamma p_{f}
\end{align}
where $p_f$ represents the fluid pressure and $\gamma$ is the fluid pressure coefficient which will be defined in a subsequent section. Note that the sign convention for the peridynamic pressure is positive in tension, while the convention for the fluid pressure is positive when it induces expansion of the bulk material.

\section{PERIDYNAMIC MATERIALS WITH EFFECTIVE STRESS}
In the present modification of the peridynamic material above, the fluid itself is not represented directly in the model, only its effect on the constitutive behavior of the porous medium. For this reason we assume that the fluid pressure is known \emph{a priori} (through alternative means of analysis for the flow quantities). It is important to recognize that the effective stress principle does not alter the constitutive behavior of the material influenced by the pore pressure. Even though the loads produced by the fluid pressure are introduced in the constitutive relationship, they can easily be extracted to appear in the formulation as an applied load. We chose to introduce the pore pressure in the constitutive relationship to simplify the numerical implementation since no method for application of natural boundary conditions has been developed yet for the peridynaimc theory. The saturated material under consideration behaves identically to its dry correlary as parameterized by its material constants, $k$ and $\mu$. For this reason, the principle of correspondence holds (the idea that a peridynamic material is an equivalent representation of a material in the classical theory). To this effect, the formulation above does not preclude application to elastic-plastic materials, even though it has been presented here for a linear material to simplify the formulation. 

The discussion above suggests that to preserve correspondence and the constitutive behavior of the peridynamic material, $\gamma$ must be defined in such a way that its origins stem from the forces induced by the interaction of the fluid with the porous material. To derive $\gamma$ in such fashion we turn to the \emph{theory of interacting continua} or \emph{mixture theory}.

\subsection{Theory of interacting continua}
The theory of interacting continua provides a method for characterizing the bulk response of a volume by incorporating the influence of the constituents of which it is composed. In the present formulation we use the theory to define an interaction force between the pore fluid and the solid matrix such that coupling takes place at the interface between them, but the effect is represented as a volumetric quantity. To present the theory of interacting continua in a consistent manner, we begin with a classical representation of the pore fluid and the porous media such that the porous media is treated as the classical equivalent of the ordinary linear elastic peridynamic material. At a later stage, we return to the peridynamic representation. 

In the theory of mixtures, the total density of the mixture and the mixture velocity are defined as:
\begin{align}
\rho = \sum_{i = 1}^N \rho^{(i)}, \qquad \boldsymbol{v} = \frac{1}{\rho} \sum_{i=1}^N \rho^{(i)} \boldsymbol{v}^{(i)}
\end{align}
where $N$ is the total number of constituents. We assume the existence of a partial traction vector, $\boldsymbol{t}^{(i)}$ and partial stress tensor, $\boldsymbol{T}^{(i)}$, such that $\boldsymbol{t}^{(i)} = \boldsymbol{T}^{(i)}\boldsymbol{n}$, where $\boldsymbol{n}$ is the normal vector to the surface. The total traction and total stress are written
\begin{align}
\boldsymbol{t} = \sum_{i = 1}^N \boldsymbol{t}^{(i)}, \qquad \boldsymbol{T} = \sum_{i=1}^N \boldsymbol{T}^{(i)} 
\end{align}
The balance of mass for each constituent is enforced via
\begin{align}
\frac{D^{(i)}\rho^{(i)}}{Dt} + \rho^{(i)}\mathrm{div}\left[ \boldsymbol{v}^{(i)}\right] = m^{(i)}
\end{align}
where $\frac{D^{(i)}(\cdot)}{Dt}$ is the total time derivative operator and the mass production terms from each constituent must balance (i.e. $\sum_{i = 1}^N m^{(i)} = 0$). The balance of linear momentum for the individual constituents is written
\begin{align}
\mathrm{div}\left[ \boldsymbol{T}^{(i)} \right] + \rho^{(i)}\boldsymbol{b}^{(i)} + \boldsymbol{c}^{(i)} = 0
\end{align}
In the above expression, $\boldsymbol{c}$ is the interaction term that represents the influence of the other constituents. Consider the case of two constituents, fluid flowing through a porous solid. We shall employ an interaction model of the following form
\begin{align}
\boldsymbol{c} = \alpha(\boldsymbol{v}^{(f)} - \boldsymbol{v}^{(s)})
\end{align}
The interaction term can be thought of as a drag term relative to the difference in the velocity of the fluid and the solid, where the interaction coefficient is $\alpha$. The model for how the fluid and solid interact can be much more complex, for example we could consider localized changes in porosity resulting form the pore pressure, but for this work we consider only the simple case of a drag-like term. To enforce a Newtonian balance of forces the interaction terms for the fluid and the solid will be of opposite sign (i.e. $\boldsymbol{c}^{(f)} = -\boldsymbol{c}^{(s)} = \boldsymbol{c}$).

For the case of a simple fluid without distortional stresses, the stress in the fluid can be written $\boldsymbol{T}^{(f)} = -p^{(f)}\boldsymbol{\mathrm{I}}.$ Note that $\boldsymbol{\mathrm{I}}$ is the identity tensor (no italics). If we assume that the velocity of the solid is much smaller than the velocity of the fluid, $\boldsymbol{v}^{(s)}$ can be neglected in the interaction term and the balance of linear momentum for the fluid becomes Darcy's law \cite{Darcy_1856}.
\begin{align}
\label{eq:fluidMomentumBalance}
\alpha \boldsymbol{v}^{(f)} + \mathrm{grad}\left[ p^{(f)}\right] = \rho^{(f)} \boldsymbol{b}^{(f)}
\end{align}
The balance of mass for the fluid can be written
\begin{align}
\frac{\partial \phi^{(s)}}{\partial t} + \mathrm{div}\left[ \boldsymbol{v}^{(f)}\right] = m^{(f)}
\end{align}

For a linear elastic solid, in the classical theory the stress tensor takes the form
\begin{align}
\boldsymbol{T}^{(s)} = \lambda^{(s)} \mathrm{div}\left[ \boldsymbol{u}^{(s)} \right] \boldsymbol{\mathrm{I}} + \mu^{(s)}\left(\mathrm{grad}\left[ \boldsymbol{u}^{(s)} \right] + \mathrm{grad}\left[ \boldsymbol{u}^{(s)} \right]^T \right)
\end{align}
where $\lambda^{(s)}$ and $\mu^{(s)}$ are the material parameters (recall that the relationship between the bulk modulus $K^{(s)}$ and $\lambda^{(s)}$ and $\mu^{(s)}$ is $K^{(s)} = \lambda^{(s)} + 2\mu^{(s)}/3$). The balance of linear momentum for the solid is written
\begin{align}
\label{eq:solidMomentumBalance}
-\alpha \boldsymbol{v}^{(f)} - \mathrm{div}\left[ \boldsymbol{T}^{(s)}\right] = \rho^{(s)} \boldsymbol{b}^{(s)}
\end{align}

In the absence of a fluid body force, equation~\eqref{eq:fluidMomentumBalance} can be solved for the fluid velocity, $\alpha \boldsymbol{v}^{(f)} = -\mathrm{grad}\left[ p^{(f)}\right]$. This result can be substituted into equation~\eqref{eq:solidMomentumBalance} to produce a modified form of the momentum balance expression for the solid
\begin{align}
\mathrm{div}\left[ \hat{\boldsymbol{T}}^{(s)}\right] + \rho^{(s)} \boldsymbol{b}^{(s)} = 0
\end{align}
where $\hat{\boldsymbol{T}}^{(s)} = \boldsymbol{T}^{(s)} - p^{(f)}\boldsymbol{\mathrm{I}}$. In the expression above we have made use of the identity $\mathrm{grad}[p^{(f)}] = \mathrm{div}[p^{(f)}\boldsymbol{\mathrm{I}}]$.

In the peridynamic theory, direct application of the above analysis implies that the fluid pressure coefficient, $\gamma$, in equation~\eqref{eq:ModifiedConstit} is 1.0, which represents exactly Terzaghi's notion of effective stress. It is well known that a unit value for the pressure coefficient in Terzaghi's method performs well for granular materials with very small ratios of the parameter, $K$, for the bulk material versus the solid constituent alone. In order to apply the proposed formulation effectively to materials like stone or marble, we introduce a pressure coefficient that represents Biot's concept of effective stress, whereby the ratio of the moduli between the bulk material and the solid constituent alone is taken into account
\begin{align}
\gamma = 1 - \frac{K}{K^{(s)}}
\end{align}
where $K$ is the bulk modulus of the volumetric material and $K^{(s)}$ is the bulk modulus of the solid constituent.
\begin{remark}
It is important to recognize that material behavior itself is unmodified. We merely chose to introduce the effect of the fluid pore pressure by means of a modified force scalar state, therefore correspondence still holds for our original ordinary linear elastic peridynamic material.
\end{remark}

\subsection{Discrete form of the governing equations}
The discrete form of the governing integral equation is written
\begin{align}
\sum_{i = 0}^{N} \left\{ \underline{\boldsymbol{T}}[\boldsymbol{x},t]\langle \boldsymbol{x}'_i - \boldsymbol{x}\rangle - \underline{\boldsymbol{T}}[\boldsymbol{x}'_i,t]\langle \boldsymbol{x} - \boldsymbol{x}'_i\rangle \right\} \Delta V_{\boldsymbol{x}'_i} + \boldsymbol{b}(\boldsymbol{x},t) = 0
\end{align}
where $N$ is the number of neighbors in the spherical neighborhood of $\boldsymbol{x}$ and $\Delta V_{\boldsymbol{x}'_i}$ is the volume of cell representing the particle at $\boldsymbol{x}'_i$. The vector force state $\underline{\boldsymbol{T}}$, is evaluated in terms of the accumulated deformation of all the particles in the neighborhood of $\boldsymbol{x}$ according to equation \eqref{eq:ForceStateDef}.

\section{NUMERICAL RESULTS}
This section presents a number of numerical examples that demonstrate the utility of the proposed method. To simplify the presentation,  for all of the examples a particular fluid pressure field is prescribed as if the flow characteristics for the given problem are known. In a forthcoming work the coupling of the proposed method with existing flow modeling frameworks is investigated in detail. Since the formulation is inherently three-dimensional, all of the one-dimensional verification problems were performed on a domain of unit thickness. The exact solutions to the verification problems were derived using partial differential equations. Although the proposed method is integral based, the results should match under certain conditions due to the principle of correspondence.

\subsection{Lighthouse consolidation}
The proposed method can be applied to consolidation of fluid saturated geologic material. The lighthouse consolidation problem represents a one-dimensional steady-state fluid-structure interaction problem with a known analytic solution. This problem is motivated by an example presented in \cite{Andersen}, and entails determining the deformation response of a partially submerged column of rock supporting a lighthouse. Including the influence of effective stress, the portion of the rock column below the water level experiences less deformation due to the pressure of the fluid in the pores. Above the water level, the rock column experiences increased subsidence due to capillary effects. For a column with cross-sectional area, $A$, length $a$ above the water level, length $b$ below the water level, supporting a lighthouse of weight $F$, the exact solution for the static vertical displacement is 
\begin{align}
    u(x) = \frac{1}{E} \left(-\psi_2 x^2\left(H(x) - \frac{1}{2}\right) - \psi_1 a x + \frac{F (b - x)}{A} + a b \psi_1 + \psi_2 \frac{ b^2}{2}\right) \qquad \forall x \in [-a,b]
\end{align}
where $\gamma_t$ is the specific weight of the rock and fluid together, $\gamma_f$ is the specific weight of the water, $E$ is the modulus of elasticity for the rock, $\psi_1$ and $\psi_2$ are given as
\begin{align}
\psi_1 = \gamma_t + \gamma_f \; \qquad
\psi_2 = \gamma_t - \gamma_f
\end{align}
and $H(x)$ is the standard heaviside function
\begin{align}
H(x) = \left\{ \begin{array}{c} 1 \qquad \mathrm{for} \; x \geq 0 \\ 0 \qquad \mathrm{for} \; x < 0 \end{array} \right\} 
\end{align}
The exact solution was obtained by integrating the one dimensional boundary value problem for axial loading of a column, using the specific weights of water and the rock as body forces. Note that this solution is the exact solution of the set of differential equations that represent classical solid mechanics. It is not precisely the exact solution to the governing equations of the peridynamic formulation above; however, the results show good performance for the proposed formulation.
\begin{figure}[htb!]
	\centering
\includegraphics[scale=0.6]{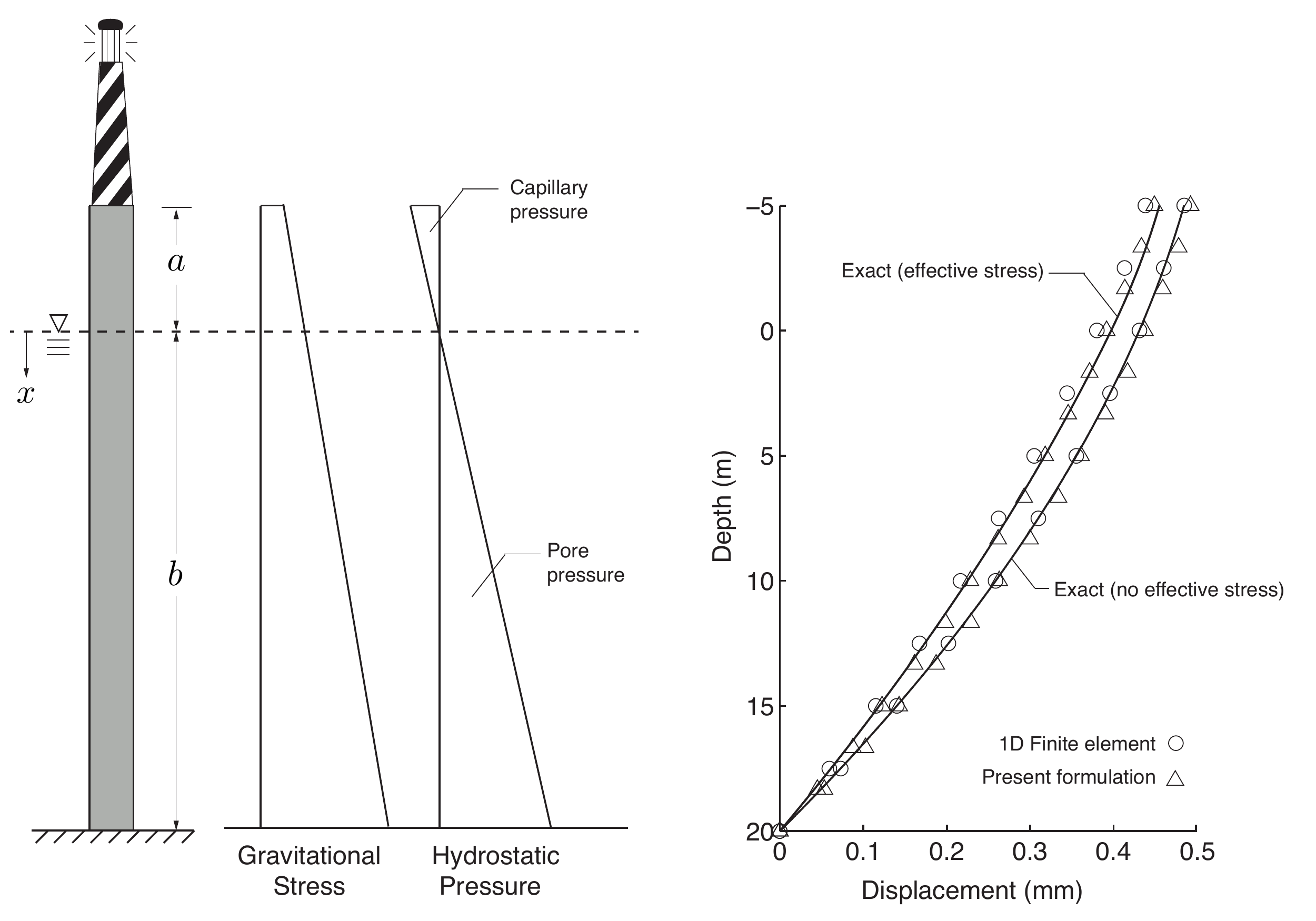}
	\caption{Lighthouse supported by a partially submerged column of rock: (left) problem geometry and applied stresses (right) verification results for the proposed formulation. The gravitational stress includes the weight of the lighthouse itself and the accumulated weight of the saturated rock column.}
	\label{fig:LighthouseDisplacements}
\end{figure}

The parameters used for this problem are given in Table~\ref{tab:LighthouseParams}. For completeness, the material constants have been presented in terms of Young's modulus, the bulk modulus, and the shear modulus. Any one of these three can be calculated from the other two.
\begin{table}
  \centering
  \begin{tabular}{cl}
    \hline
    Parameter & Value\\
    \hline
    $a$ & 5 m\\
    $b$ & 20 m \\
    $F$ & 20 MN \\
    $A$ & 78.54 $\mathrm{m}^2$  \\
    $\gamma_t$ & 25.67 $\mathrm{kN/m^3}$\\
    $\gamma_f$ & 10.02 $\mathrm{kN/m^3}$\\
    $E$ & 29.7 GPa \\
    $k$ & 9.0 GPa \\
    $\mu$ & 15.0 GPa \\
    $\gamma$ & 1.0 \\
    $\delta$ & 5.0 \\
    \hline \\
    \\
  \end{tabular}
  \caption{Lighthouse consolidation problem parameters}
  \label{tab:LighthouseParams}
\end{table}
The peridynamic analysis was performed for a 25 m long rectangular column of unit cross-sectional area discretized with elements of a volume 0.0015625 $\mathrm{m^3}$. The boundary conditions consist of fixing the bottom layer of elements such that the vertical displacement is zero and applying a concentrated load to each element in the top layer consistent with the weight of the lighthouse distributed over the number of elements in the top layer. The gravitational force was applied via the body force term $\boldsymbol{b}(x)$ and the fluid pressure was applied via the $p_f$ term in equation~\eqref{eq:ModifiedConstit}. Figure~\ref{fig:LighthouseDisplacements} shows a comparison of the results for this problem using the proposed formulation and linear, one-dimensional, two-node bar finite elements. Although the error for both methods is relatively small, the finite element solution begins to deviate from the exact solution for the case that includes effective stress.

The results illustrate that the modification of the pressure term according to equation~\eqref{eq:ModifiedConstit} does not invalidate the principle of correspondence. Even though in the proposed formulation the influence of the fluid pore pressure is engendered though the calculation of the force scalar state, the underlying material responds to loading in the same way as described by classical constitutive theory.

\subsection{Time dependent consolidation due to harmonic loading}
\begin{figure}[htb!]
	\centering
\includegraphics[scale=0.65]{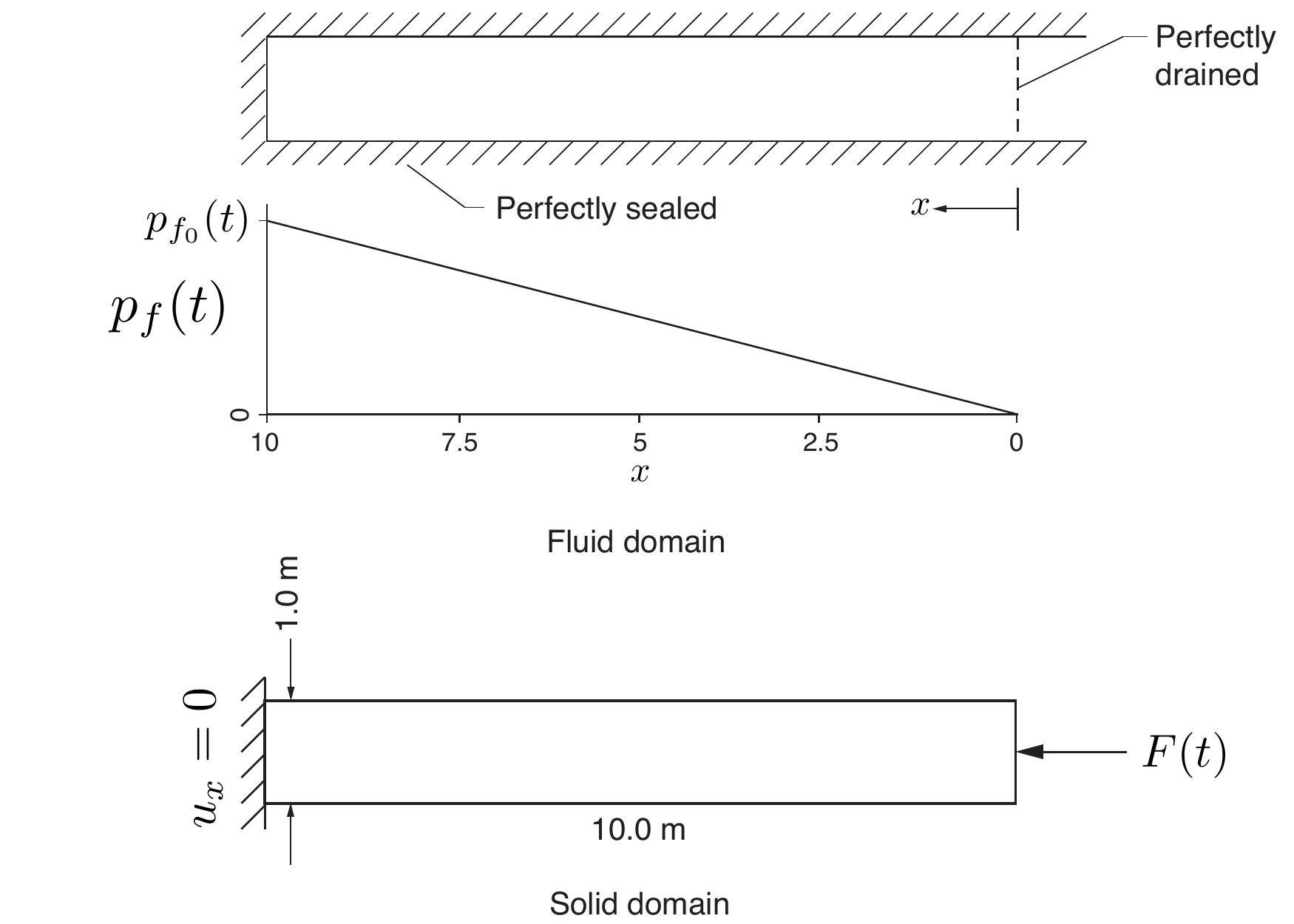}
	\caption{Problem geometry, boundary conditions and loading for the time dependent consolidation problem.}
	\label{fig:TimeDepConDomain}
\end{figure}

For additional verification, we consider here the application of the proposed method to time dependent consolidation of a fluid saturated porous material under harmonic loading. This one dimensional, shown in Figure~\ref{fig:TimeDepConDomain}, consists of a column of material fixed at one end and loaded with an excitation force, $F(t)$, at the other given by
\begin{align}
F(t) = 50(1 - \mathrm{cos}(75t))
\end{align}
The fluid pressure is assumed to vary according to the excitation force in time
\begin{align}
p_f(t) = 2F(t)\left(1 - \frac{t}{\tau}\right)
\end{align}
and linearly along the domain dropping from $p_f$ at the left side to zero at the right side where the excitation force is applied. These conditions correspond to a perfectly drained boundary at the surface of loading. At time $\tau$ the pressure field reaches a state of zero pressure. The peridynamic analysis for this problem was performed on a unit-cross sectional domain of depth, $L=10$ m discretized into cubic cells of volume 0.0026653 $\mathrm{m}^3$. The problem parameters are listed in Table~\ref{tab:TimeDepConParams}.
The exact solution for the deflection of the surface on which the loading is applied can be obtained easily by integrating the one dimensional boundary value problem of classical elasticity and is given as
\begin{align}
u_{x=0}(t) = \frac{L}{E}\left(F(t) + \frac{p_{f_0}(t)}{2}\right)
\end{align}

\begin{table}
  \centering
  \begin{tabular}{cl}
    \hline
    Parameter & Value\\
    \hline
    $F$ & $50(1 - \mathrm{cos}(75 t) ) $kN \\
    $\tau$ & 0.4 s \\
    $E$ & 10.0 MPa \\
    $k$ & 3.33 MPa \\
    $\mu$ & 5.0 MPa \\
    $\gamma$ & 1.0 \\
    $\delta$ & 3.5 \\
    \hline \\
    \\
  \end{tabular}
  \caption{Time dependent consolidation due to harmonic loading problem parameters}
  \label{tab:TimeDepConParams}
\end{table}

\begin{figure}[htb!]
	\centering
\includegraphics[scale=0.6]{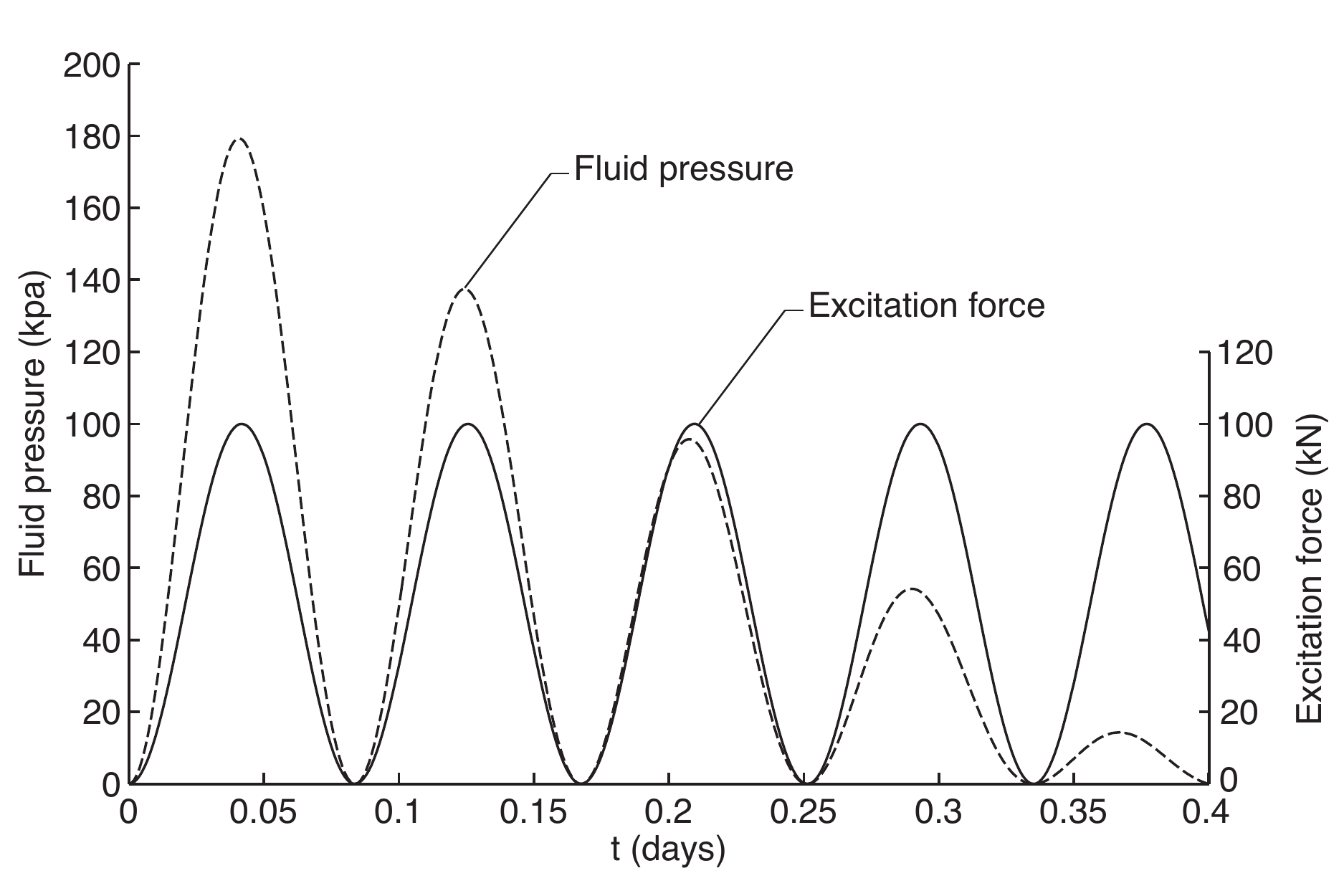}
	\caption{Evolution of the excitation force and fluid pressure in time for the time dependent consolidation due to harmonic loading problem. The fluid pressure is measured at the maximum depth of the domain ($x = 10$ m) and the excitation force is measured at the top surface of the domain ($x = 0$ m).}
	\label{fig:TimeDepConForces}
\end{figure}

\begin{figure}[htb!]
	\centering
\includegraphics[scale=0.6]{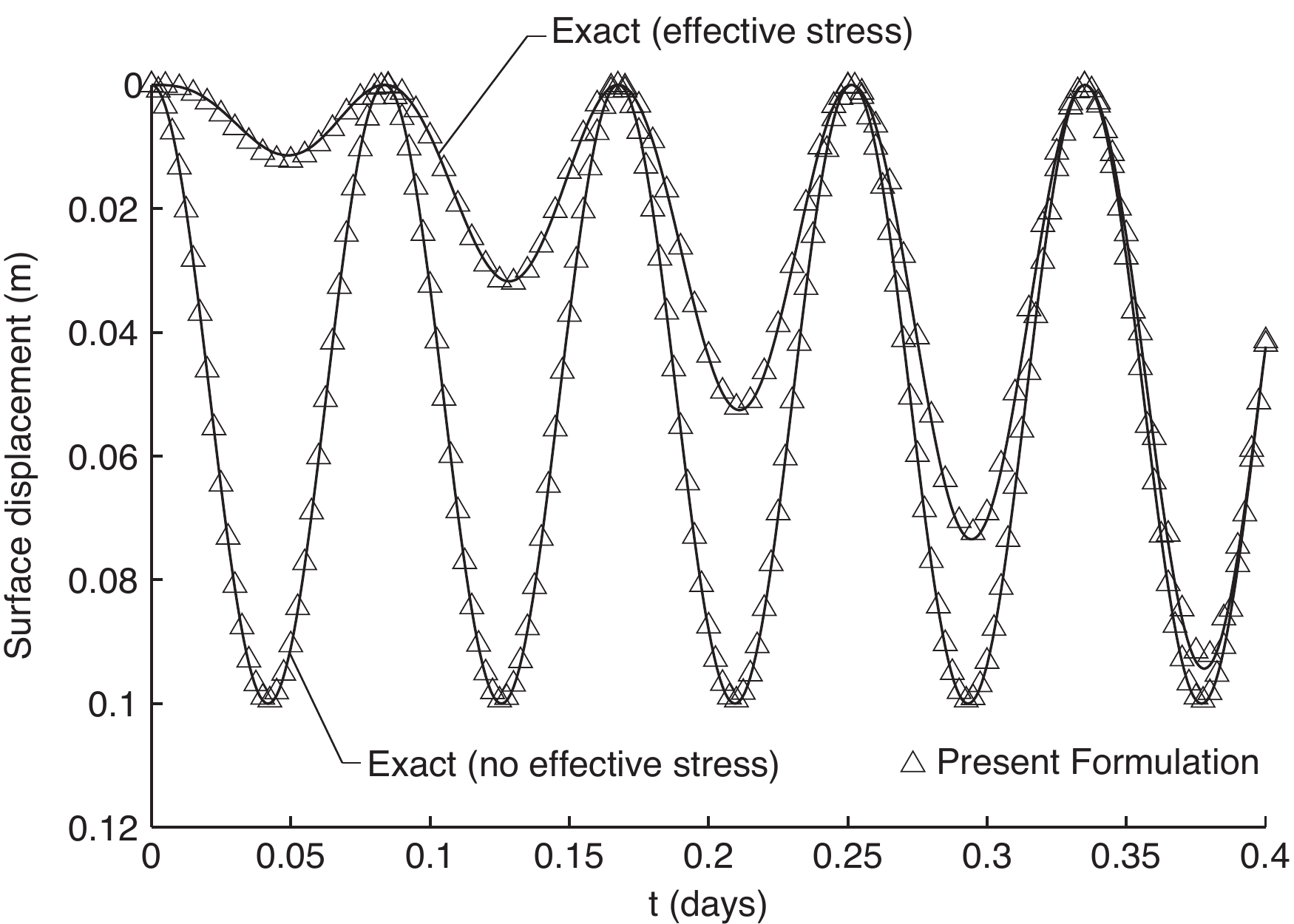}
	\caption{Evolution of surface deflection for the time dependent consolidation due to harmonic loading problem comparing the behavior including the effect of pore pressure and the behavior of a dry material.}
	\label{fig:TimeDepCon}
\end{figure}

Figure~\ref{fig:TimeDepConForces} shows the harmonic loading and the maximum value of the fluid pressure, each measured at opposite ends of the domain. Figure~\ref{fig:TimeDepCon} shows the evolution of surface deflection in time for both a dry porous material and one in which effective stress is taken into account. For both cases the proposed method performs well in capturing time dependent consolidation.


\subsection{Surface subsidence due to well extraction}
The next problem demonstrates the surface deflection caused by a decrease in fluid pressure inside a reservoir as fluid is extracted. The drop in pore pressure from an initial state of equilibrium causes the porous material to contract volumetrically leading to subsidence. This problem is motivated by a similar problem studied in \cite{Dean, MartinezLDRD}. The fluid pressure field for this problem is based on the numerical results of these studies. A linear drop in pressure is assumed form the borehole radially to the outer boundary of the domain. The change in pressure from extraction is assumed to be -150 psi at the borehole which represents an extraction duration of 75 days. The authors of \cite{Dean, MartinezLDRD} report that the surface subsidence should be 4.5 ft after 75 days.

\begin{figure}[htb!]
	\centering
\includegraphics[scale=0.8]{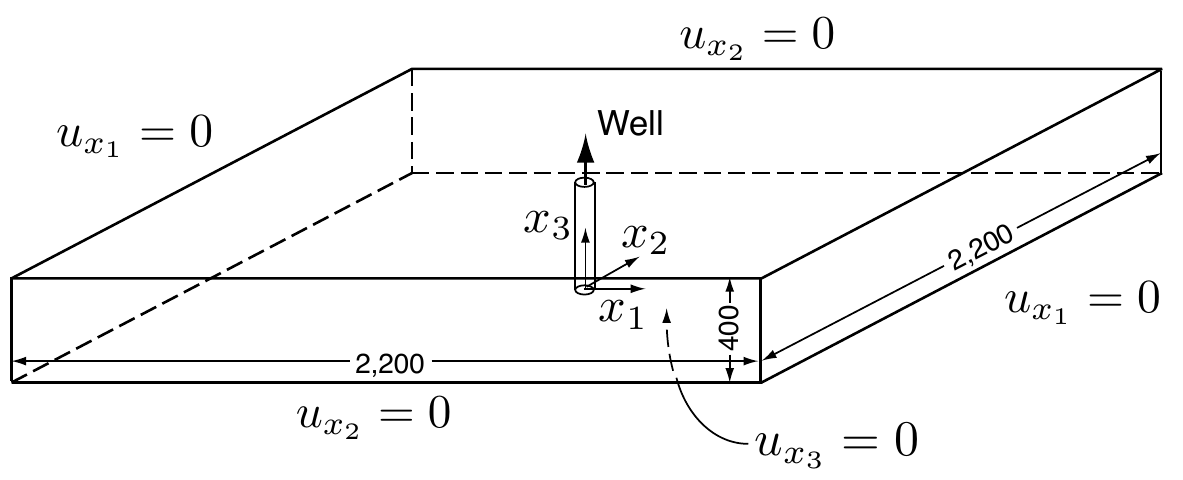}
	\caption{Problem geometry and boundary conditions for the surface subsidence due to well extraction problem (all dimensions are in ft).}
	\label{fig:SufaceSubsidenceDomain}
\end{figure}

\begin{figure}[htb!]
	\centering
\includegraphics[scale=0.6]{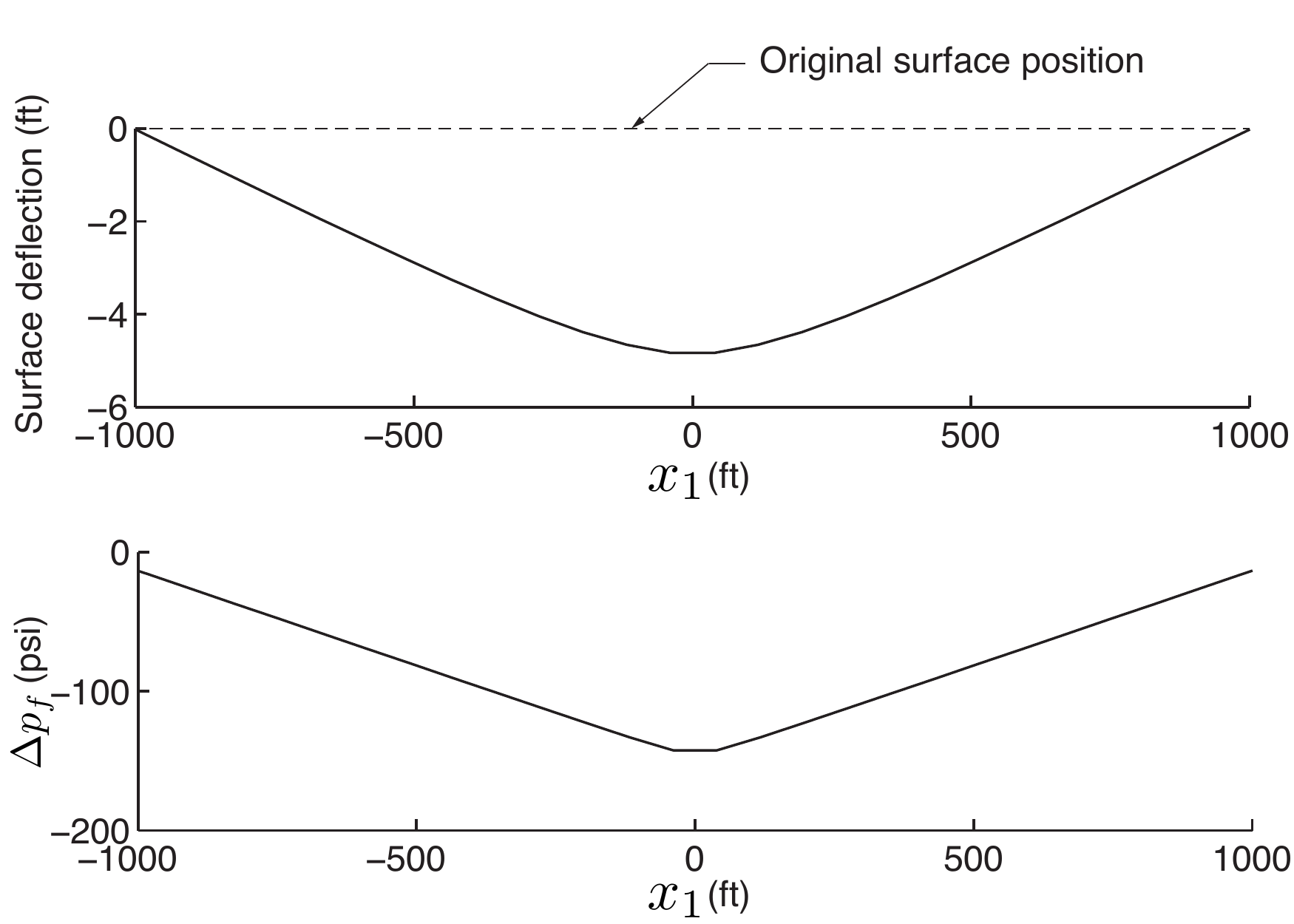}
	\caption{Surface subsidence due to well extraction. Extraction from the reservoir leads to a decrease in fluid pressure near the borehole. This in turn results in consolidation of the reservoir material. The displacement in the $x_3$-direction is measured along the $x_3$-axis at the top surface of the domain. The loading conditions correspond to an extraction duration of 75 days.}
	\label{fig:SurfaceSubsidence}
\end{figure}

The peridynamic analysis of this problem using the proposed method was performed on a domain of dimensions shown in Figure~\ref{fig:SufaceSubsidenceDomain}. The field-scale domain was discretized with cells of inner radius of 40 ft and the parameters used are reported in Table~\ref{tab:SurfaceSubsidenceParams}. The boundary conditions consist of symmetry conditions for displacement along the outer edges of the domain and zero displacement in the $x_3$ direction along the bottom surface. From the results shown in Figure~\ref{fig:SurfaceSubsidence} it is clear that the proposed formulation captures the same effective subsidence of the surface due to extraction as reported in previous studies.

\begin{table}
  \centering
  \begin{tabular}{cl}
    \hline
    Parameter & Value\\
    \hline
    $k$ & 3.8642 ksi \\
    $\mu$ & 3.846.2 ksi \\
    $\Delta p$ & -150.0 psi \\
    $\gamma$ & 1.0 \\
    $\delta$ & 3.5 \\
    \hline \\
    \\
  \end{tabular}
  \caption{Surface subsidence due to well extraction problem parameters}
  \label{tab:SurfaceSubsidenceParams}
\end{table}

\subsection{The effect of leak-off during hydraulic fracturing}
With regard to applications in hydraulic fracturing, the present formulation is well suited for capturing the influence of leak-off effects. As fracturing fluid is injected into the reservoir, depending on the properties of the fluid and the injection rate, the fluid will penetrate into the geologic material to a varying extent. For high injection rates or geologic materials with low porosity very little leak-off occurs and the primary driver for damage propagation comes directly from the pressure of the fluid on the surface of the borehole. When the injection rate is lower or the penetrability of the material itself is high the fluid leaks into the pores of the geologic material leading to an increased state of stress in the porous media. To avoid the complexity of the interaction between crack growth and fluid penetration we only consider here the process for pristine material prior to crack initiation. The effective influence of leak-off on the porous material is measured by the dilation of borehole in terms of displacement for various assumed depths of fluid penetration.

\begin{figure}[htb!]
	\centering
\includegraphics[scale=0.6]{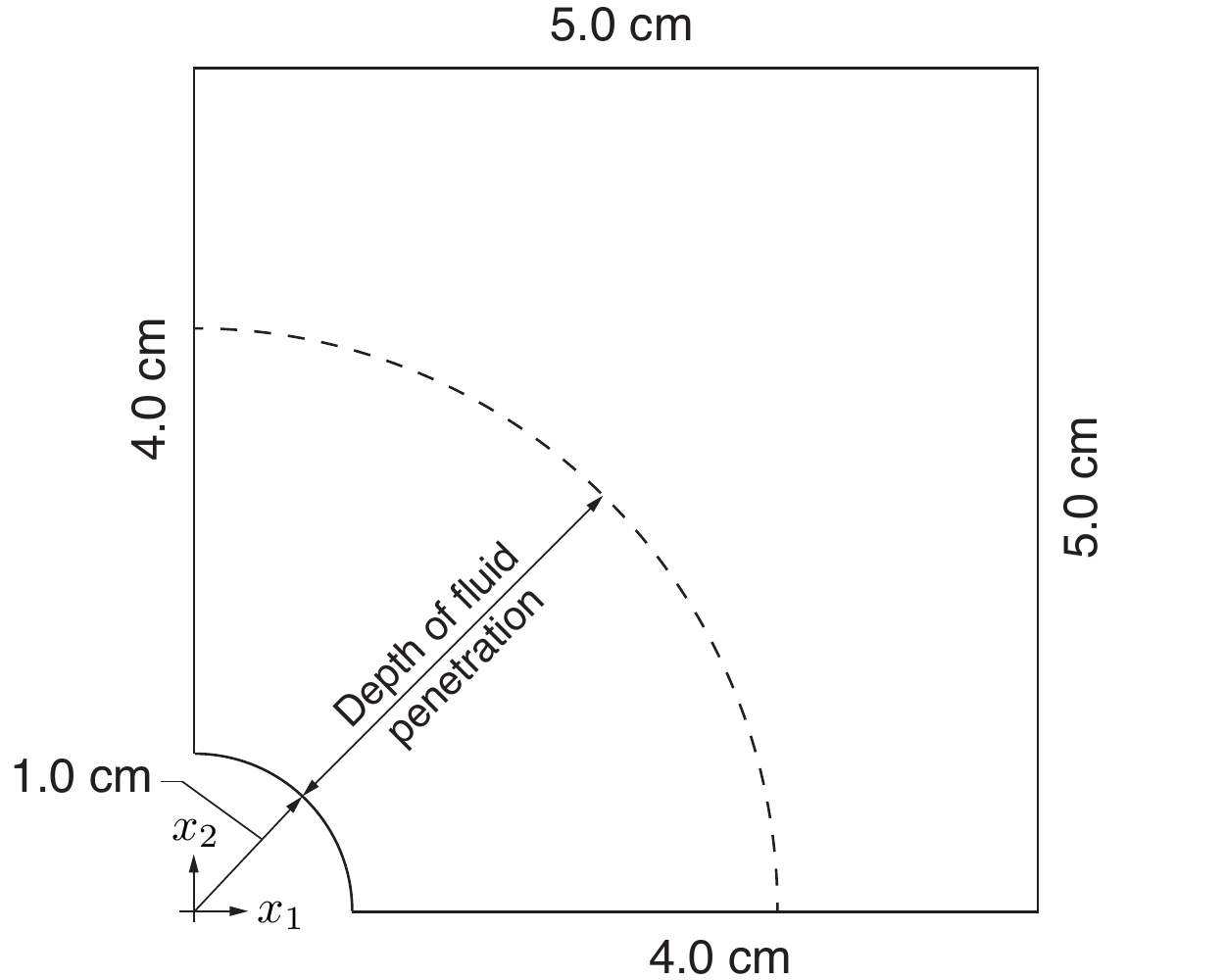}
	\caption{Geometry of the leak-off problem. In the region of fluid penetration, the pressure is specified as a constant value at the borehole surface and zero at the outer boundary. Between these two surfaces, a linear drop in pressure is prescribed. There is no gradient in the pressure in the $x_3$ direction.}
	\label{fig:PlateHoleDiagram}
\end{figure}

The geometry of the leak-off problem is shown in Figure~\ref{fig:PlateHoleDiagram}. The domain was discretized with cells of volume 0.001153 $\mathrm{cm^3}$. The analysis was performed on a three dimensional domain with a thickness of 0.5 cm assuming plane-strain conditions (no strain in or out of the plane of the paper). Symmetry was used to reduce the problem to model only one quarter of the material in the vicinity of the borehole. The loading was prescribed as a constant force applied to each element equivalent to 1.0 Mpa of confining pressure on the top and right surface. On the borehole surface a force was prescribed on each element equivalent to a pressure of 0.1 Gpa. To preserve symmetry, no displacement was allowed in the $x_2$-direction along the bottom surface and no displacement was allowed in the $x_1$-direction along the left surface. The displacement was prescribed as zero in the $x_3$-direction on the top and bottom surface.

\begin{table}
  \centering
  \begin{tabular}{cl}
    \hline
    Parameter & Value\\
    \hline
    $k$ & 9.0 GPa \\
    $\mu$ & 15.0 GPa  \\
    $\gamma$ & 1.0 \\
    $\delta$ & 3.5 \\
    \hline \\
    \\
  \end{tabular}
  \caption{Hydraulic fracturing problem parameters}
  \label{tab:HydraulicFracturingParams}
\end{table}

\begin{figure}[htb!]
	\centering
\includegraphics[scale=0.65]{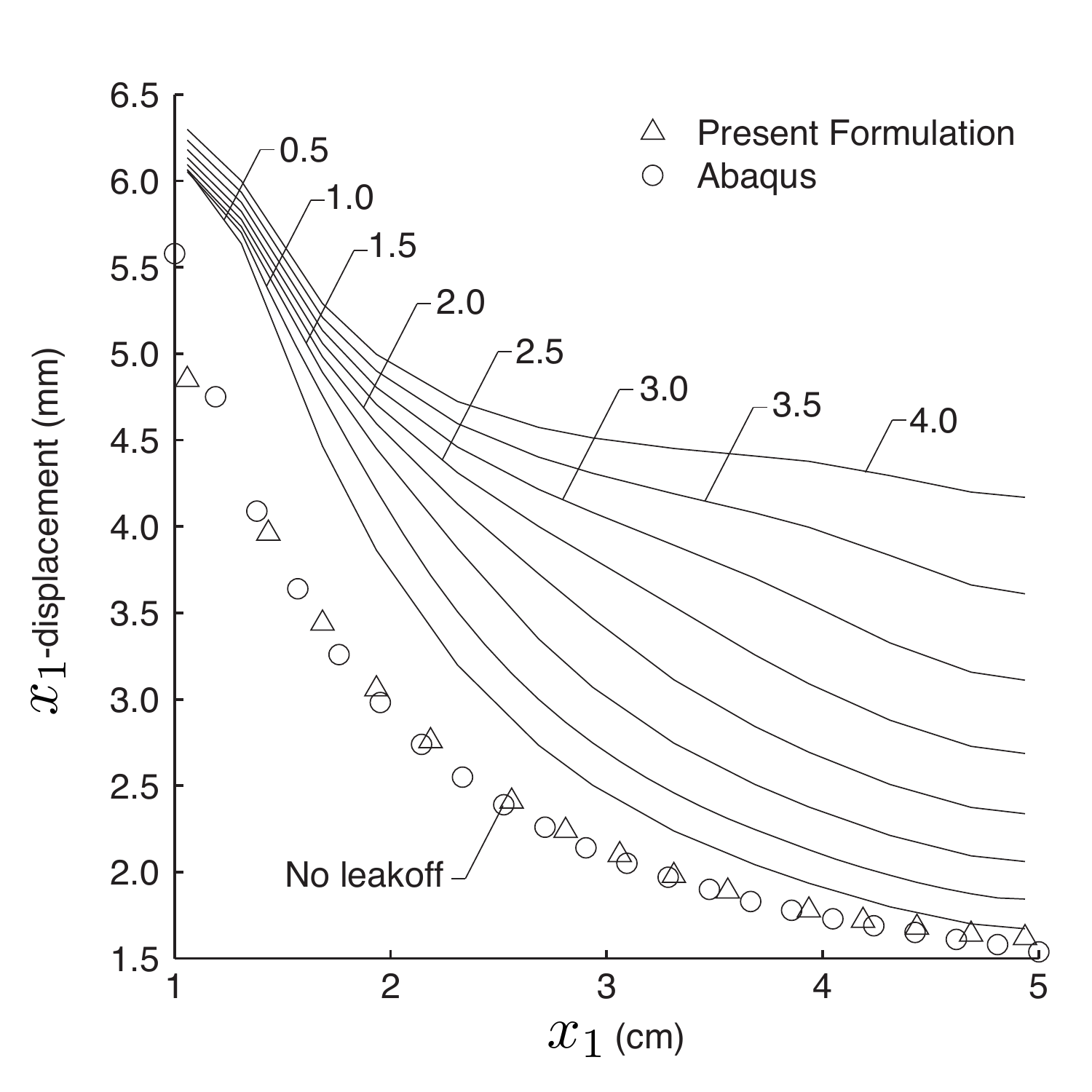}
	\caption{The effect of leak-off on the displacement in the vicinity of the borehole during hydraulic fracturing. The symbols show verification of the case without leak-off by comparison with Abaqus. The solid lines represent the effect of leak-off with varying depths (in centimeters) of fluid penetration beyond the borehole surface.}
	\label{fig:PlateHoleResults}
\end{figure}

To investigate the effect of leak-off several simulations were run with varying degrees of fluid penetration. For each case, the fluid pressure was prescribed as 0.1 Gpa on the borehole surface and zero at the depth of penetration with a linear drop in pressure between. The results of this study are shown in Figure~\ref{fig:PlateHoleResults} which plots the displacement in the $x_1$-direction along the $x_1$-axis. To verify the solution comparisons are made, for the case of no fluid penetration with the results produced by Abaqus (a commercial finite element code) using eight-node brick elements with similar boundary conditions, loading, and element density. Good correspondence is obtained for this case between the two codes.

Also shown in Figure~\ref{fig:PlateHoleResults} is the effect of varying the depth of fluid penetration. Near the borehole, the displacement is greatly increased due to leak-off effects for slight penetration, but the influence does not continue to grow for greater penetration depths. However, the response away from the borehole is almost unaffected by the initial penetration, but continues to increase rapidly as the fluid depth increases. Representative values of the displacement both at the borehole surface and the right edge of the domain are given in Table~\ref{tab:LeakoffDispValues} for the dry case and the case of 4.0 cm of fluid penetration.

\begin{table}
  \centering
  \begin{tabular}{ccc}
    \hline
    & Borehole Edge & Right Domain Boundary\\
    & ($x_1 = 1.0$) & ($x_1=5.0$) \\
    \hline
    Dry & 4.85 mm & 1.62 mm \\
    4.0 cm of leakoff & 6.30 mm & 4.17 mm\\
    \hline \\
    \\
  \end{tabular}
  \caption{Comparison of the displacements for the dry case and the case with 4.0 cm of fluid penetration}
  \label{tab:LeakoffDispValues}
\end{table}

If the pore fluid pressure is large enough, in comparison to the confinement pressure due to in-situ conditions, the presence of the fluid leads to increased stress in the porous material rather than stress relief (as was seen in the lighthouse consolidation example) and becomes the driving force behind damage propagation. Cases like this represent a fundamental deviation from the context of consolidation or subsidence problems for which the traditional effective stress principle was developed.

\section{SUMMARY}
This work presents a nonlocal model for fluid-structure interaction based on an effective stress concept. Through the principle of correspondence the underlying constitutive behavior of the porous material is unaltered, but the effect of fluid pore pressure is accurately captured. The method was demonstrated for a number of numerical examples including, one-dimensional consolidation and three-dimensional subsidence. the proposed formulation was also used to show the effect of leak-off during the pre-crack initiation stages of hydraulic fracture. The results of each of these examples show strong correlation with results previously reported and suggest that the formulation accurately captures the effective stress principle. This outcome is significant from the standpoint that it enables the development of more advanced methods for fluid-structure interaction that incorporate all the advantages of a non-local formulation such as discontinuity capturing and long range forces. Ultimately, it is intended that this work lay the foundation for the study of more complex interaction caused by flow induced damage of porous materials which can be used to model engineering problems like hydraulic fracturing.

\bibliographystyle{unsrt}
\bibliography{../Master_References}

\end{document}